# Enhancing Optimization Performance: A Novel Hybridization of Gaussian Crunching Search and Powell's Method for Derivative-Free Optimization

Benny Wong

Techdesignism Machine Learning Laboratory

Abstract: This research paper presents a novel approach to enhance optimization performance through the hybridization of Gaussian Crunching Search (GCS) and Powell's Method for derivative-free optimization. While GCS has shown promise in overcoming challenges faced by traditional derivative-free optimization methods [1], it may not always excel in finding the local minimum. On the other hand, some traditional methods may have better performance in this regard. However, GCS demonstrates its strength in escaping the trap of local minima and approaching the global minima. Through experimentation, we discovered that by combining GCS with certain traditional derivative-free optimization methods, we can significantly boost performance while retaining the respective advantages of each method. This hybrid approach opens up new possibilities for optimizing complex systems and finding optimal solutions in a range of applications.

## Introduction

Optimization is a fundamental task in various fields, aiming to find the best possible solution for a given problem. Traditional optimization methods often rely on derivative-based approaches, which require the availability of analytical derivatives. However, in many real-world scenarios, obtaining these derivatives may be challenging or even impossible. As a result, derivative-free optimization methods have gained significant attention due to their ability to tackle such scenarios.

Among derivative-free optimization methods, Gaussian Crunching Search (GCS) has emerged as a promising approach. GCS utilizes a probabilistic model to guide the search process, making it well-suited for problems with complex and non-linear landscapes. It has shown effectiveness in overcoming some of the difficulties faced by traditional derivative-free optimization methods.

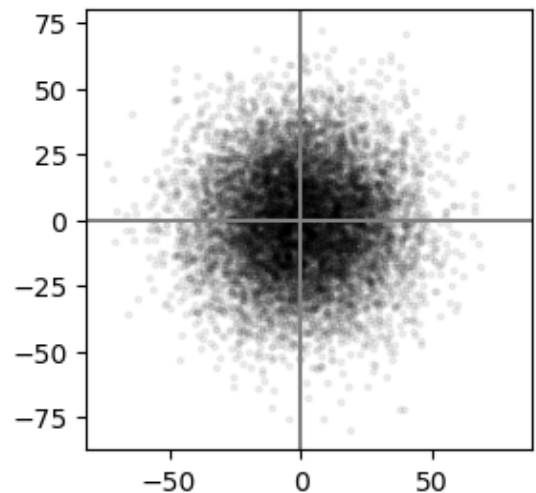

Figure 1: GCS Stochastic Mutation [1]

However, while GCS has demonstrated its strengths in escaping the trap of local minima and approaching global minima, it may not always excel in finding the local minimum. On the other hand, certain traditional derivative-free optimization methods have shown better performance in this aspect. Therefore, there is a need to explore the potential benefits of combining GCS with these traditional methods to enhance optimization performance.

In this research paper, we propose a novel hybridization approach that combines GCS with some traditional derivative-free optimization methods, specifically Powell's Method[2]. By blending the strengths of GCS and Powell's Method, we aim to improve optimization performance in terms of finding both local and global minima. The hybridization approach allows us to leverage the advantages of each method, resulting in a more robust and efficient optimization process.

Through extensive experimentation and evaluation, we investigate the performance of the hybridized GCS and Powell's Method in various optimization tasks. We compare the hybrid approach against standalone GCS and traditional derivative-free optimization methods to assess its effectiveness in terms of convergence speed and solution quality. Our findings reveal that the hybridization of GCS and traditional methods can significantly boost optimization performance while retaining their respective advantages.

## Local minima and High-dimension

Traditional derivative-free optimization methods often encounter difficulties when dealing with local minima, which can impede their ability to find the global optimum. These methods aim to identify the optimal solution by iteratively exploring the search space without relying on derivative information. However, in the presence of multiple local minima, traditional derivative-free methods may become trapped in suboptimal solutions, preventing them from reaching the global minimum. Furthermore, as the dimensionality of the optimization problem increases, the search space becomes exponentially larger, leading to the curse of dimensionality. This phenomenon can cause traditional methods to struggle in high-dimensional optimization tasks, resulting in poor performance and slow convergence rates.

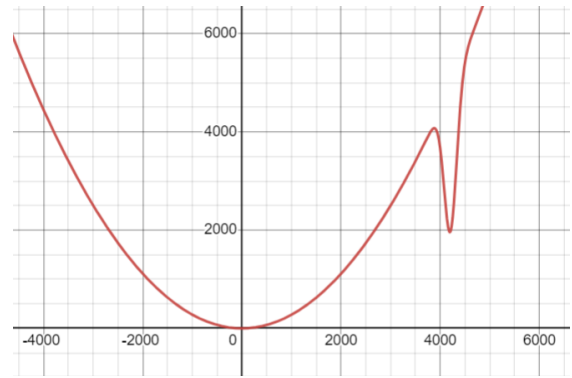

Figure 2: Global Min. & Local Min.

Gaussian Crunching Search (GCS) offers a promising solution to these challenges. GCS is a novel derivative-free optimization method that utilizes a probabilistic model to guide the search process. By incorporating randomness and adaptability, GCS can effectively escape local minima and explore the search space more efficiently. Unlike traditional methods, GCS does not rely on derivative information, making it particularly suitable for high-dimensional optimization problems where obtaining such information is often difficult or impractical. The probabilistic nature of GCS allows it to dynamically adjust its search strategy, enabling it to explore promising regions and avoid becoming trapped in local minima. Through its unique approach, GCS has demonstrated the capability to overcome the limitations of traditional derivative-free optimization methods, making it a valuable tool for addressing the challenges associated with local minima and high-dimensional optimization.

## Benchmark

When considering benchmark problems, there are three crucial aspects that need to be taken into account: the presence of a single local minimum, the existence of multiple local minima, and the high-dimensional nature of the problem[3].

### Benchmark 1

Firstly, benchmark problems with a single local minimum serve as a standard test to assess an optimization algorithm's ability to converge towards the global optimum. These problems are characterized by having a unique solution that optimizes the objective function. Evaluating the performance of an algorithm on such problems provides insights into its ability to efficiently locate the global optimum and avoid getting stuck in local minima.

$$f_1(x) = \left(1 + 19 - 19e^{-\frac{\left(-20 + \frac{x}{60}\right)^2}{9}}\right)\left(\frac{x}{60}\right)^2$$

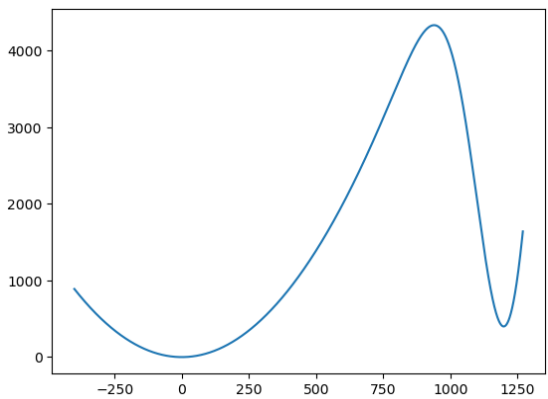

Figure 3: Benchmark function 1

Results:

| Method | Optima [x,y] |
|--------|--------------|
| P-GCS  | [1.22e-33, 8.314e-69] |
| GCS    | [0.44, 0.00108] |
| Powell | [1198.58, 399.53] |

Both the Powell-Gaussian Crunching Search (P-GCS) and Gaussian Crunching Search (GCS) algorithms exhibit the ability to escape local minima and approach the global minimum located at coordinates (0,0). However, a notable distinction arises between these two algorithms in terms of their effectiveness in achieving optimal solutions. While both P-GCS and GCS can successfully escape local minima, P-GCS demonstrates a superior capability to approach a better optimum that is closer to the global minimum. On the other hand, Powell's method tends to become trapped at local minima, limiting its ability to converge towards the global minimum. This observation emphasizes the advantage of incorporating GCS into the hybridization framework, as it allows for better optimization performance and an increased likelihood of reaching an optimal solution. By combining the strengths of Powell and GCS, the hybrid approach(P-GCS) provides a promising solution for addressing the challenges associated with local minima and achieving more favourable outcomes in derivative-free optimization problems.

**Benchmark 2**

Secondly, benchmark problems with multiple local minima are designed to challenge optimization algorithms by introducing situations where the global optimum is not easily reachable. These problems often exhibit complex landscapes with numerous local minima, making it difficult for algorithms to identify the global optimum. Evaluating the performance of an algorithm on such problems helps assess its ability to explore the search space effectively and escape from local minima to find the global optimum.

$$f_2(x) = -\lambda e^{-\mu\sqrt{x^2+y^2}} + \lambda - cos(x) - cos(y) + 2$$

$\lambda = 15, \mu = 0.05$

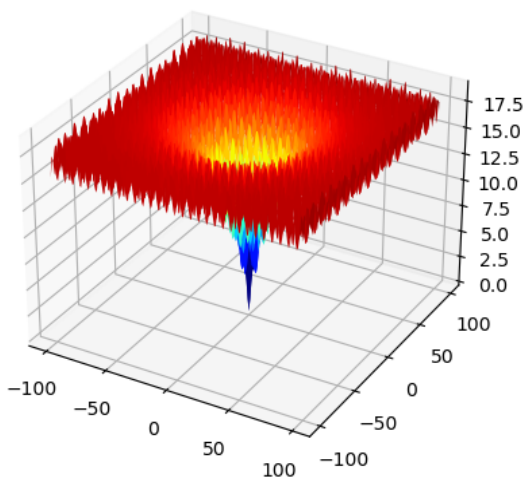

Figure 4: Benchmark function 2

Results:

| Method | Optima [x,y,z] |
|--------|----------------|
| P-GCS  | [1.97e-12, -1.21e-12, 1.73e-12] |
| GCS    | [-0.01, -0.02, 0.0143] |
| Powell | [603.19, 603.19, 15.] |

In benchmark test 2, where there are multiple local minima, the optimization task becomes more challenging compared to test 1.

However, both the Powell-Gaussian Crunching Search (P-GCS) and Gaussian Crunching Search (GCS) algorithms exhibit the ability to escape local minima and approach the global minimum located at coordinates (0,0). It is noteworthy that P-GCS demonstrates particular strength in approaching a better optimum that is closer to the global minimum, indicating its effectiveness in navigating complex landscapes with multiple local minima. In contrast, Powell's method tends to become trapped at local minima, limiting its ability to converge towards the global minimum. This observation highlights the advantage of incorporating GCS into the hybridization framework, as it enables better optimization performance and increases the likelihood of achieving more favourable outcomes in scenarios with multiple local minima. By combining the strengths of Powell and GCS, the hybrid approach offers a promising solution for addressing the challenges associated with optimizing problems that involve multiple local minima.

**Benchmark 3**

Lastly, benchmark problems that involve high-dimensional spaces pose a different set of challenges for optimization algorithms. As the dimensionality of the problem increases, the search space becomes exponentially larger, making it more challenging for algorithms to explore and find the optimal solution. Evaluating the performance of an algorithm on high-dimensional benchmark problems allows for an assessment of its scalability and efficiency in handling real-world problems with a large number of variables.

$$f_3(x_0, x_1, \ldots, x_k) = -\lambda e^{-\mu \sqrt{\sum_{i=0}^{k} x_i^2}} + \lambda$$

$$\lambda = 15, \mu = 0.05$$

Figure 5: Benchmark function 3

Results:

| Method | Optima |
|--------|--------|
| P-GCS | [-8.60e-13, 2.02e-12, 3.52e-12, -5.47e-12, -4.60e-12, -2.05e-12, -2.12e-12, 1.11e-12, -1.88e-12, -2.68e-12, 4.98e-12 2.01e-12, **8.13e-12**] |
| GCS | [-0.03, -0., 0., 0.01, -0.01, -0.01, 0.02, -0.01, -0.02, 0.02, 0., -0.01, **0.035**] |
| Powell | [200.63, 198.4, 201.64, 198.4, 201.64, 198.4, 201.02, 198.4, 202.5, 198.4, 201.02, 198.4, **15.**] |

In benchmark test 3, the focus is on evaluating the performance of optimization methods in high-dimensional spaces. This scenario poses a significant challenge for algorithms due to the exponential increase in the search space as the dimensionality increases. However, when considering the performance of the Powell-Gaussian Crunching Search (P-GCS) and Gaussian Crunching Search (GCS) algorithms, it becomes evident that P-GCS excels in approaching a better optima that is closer to the global minimum. P-GCS demonstrates its effectiveness in navigating complex, high-dimensional landscapes by efficiently exploring the search space and converging towards more favourable solutions. On the other hand, Powell's method tends to become trapped at local minima, limiting its ability to escape and approach the global minimum in high-dimensional spaces. This observation further emphasizes the advantage of incorporating GCS into the hybridization framework, as it helps overcome the limitations of Powell's method and enhances the optimization performance, particularly in high-dimensional scenarios.

Considering these three aspects of benchmark problems – single local minimum, multiple local minima, and high-dimensionality – provides a comprehensive evaluation of an optimization algorithm's performance. By assessing an algorithm's ability to handle these varying scenarios, researchers can gain insights into its strengths and limitations, enabling further improvements and advancements in the field of optimization.

## The new equation in P-GCS

In our efforts to enhance the crunching process, we have developed a new equation called the Bounded Wave equation (BW equation). This innovative equation allows for more efficient control over the periodic variation. By incorporating the Bounded Wave equation into the crunching process, we are able to better manage and manipulate the periodicity of the variation, resulting in improved optimization performance. This advancement enables us to fine-tune the search process and achieve more accurate and effective results. The Bounded Wave equation represents a significant improvement in optimizing algorithms, providing a valuable tool for researchers and practitioners in various fields.

**Bounded Wave equation (BW equation):**

$$w(x,a,b) = t(s(x,a,b))$$

where

$$t(x) = \tan\left(\pi x - \frac{\pi}{2}\right)^2$$

$$T(x) = \frac{\text{arccot}(\sqrt{x})}{\pi}$$

$$s(x,a,b) = \frac{\cos(2\pi x - \pi) + 1}{2}(T(b) - T(a)) + T(a)$$

In the Bounded Wave equation, the parameters 'a' and 'b' are instrumental in determining the bounds of the function. These parameters allow us to define the range within which the equation operates, providing control over the behaviour of the equation. Additionally, the output value of the Bounded Wave equation serves a crucial role in controlling the standard deviation of the mutation in the Gaussian Crunching Search (GCS) algorithm. By utilizing the output value as a parameter for mutation, we can modulate the extent of exploration and exploitation in the search process. Notably, the parameter 'b' in the Bounded Wave equation can be set to infinity, which proves to be particularly useful in finding better optima and escaping local minima. By allowing the function to have an unbounded upper limit, the algorithm can explore a wider range of solutions and increase the likelihood of finding more favourable optima. This flexibility in setting the bounds of the function, along with the ability to control the standard deviation of mutation, empowers the optimization process and improves the overall performance of GCS.

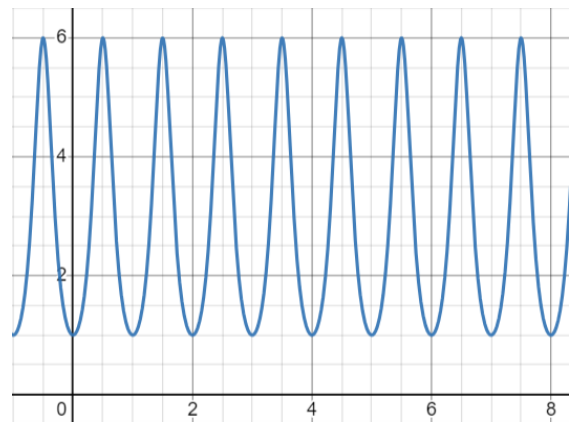

$w(x, 1, 6)$

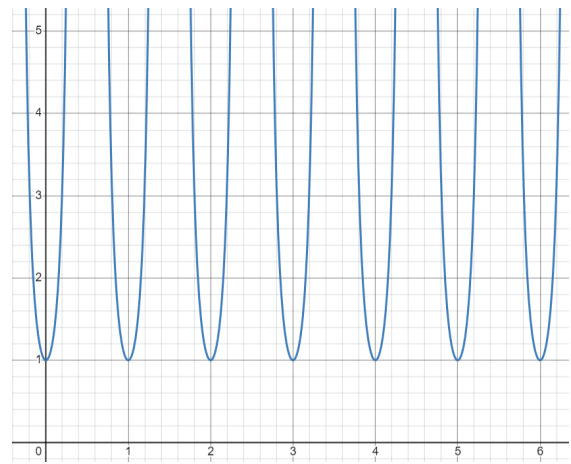

$w(x, 1, \infty)$

# Implementation

Python Code (P-GCS):

```python
#P-GCS
from numpy import *
from scipy.optimize import minimize
import scipy
import numpy as np
np.set_printoptions(precision=2)
import datetime

def objective(x_list):
    ld = 15
    mu = 0.05

    return -ld* np.exp(-mu* np.sqrt( \
       np.sum( np.square(np.array(x_list)) ) ))  + ld

def t(x):
    return np.square( np.tan(np.pi*x-np.pi/2) )

def T(x):
    return (np.pi/2-np.arctan( np.sqrt(x) ) )/np.pi

def s(x,a,b):
    return ((np.cos(2*np.pi*x-np.pi)+1)/2)*(T(b)-T(a))+T(a)

def w(x,a,b):
    return t(s(x,a,b))

target = 0.05

start_point = [200, 200, 200, 200, 200, 200, 200, 200, 200, 200, 200, 200]

##############################################
start_time = ( datetime.datetime.now() )

w_a = 0
w_b = np.inf

sd = w_a
current_pos = start_point
current_obj_value = objective(current_pos)

crunch_period = 5000
crunch_step = 0

w_cache = np.array([])

for i in range(crunch_period):
    w_cache = np.append(w_cache, w(i/crunch_period, w_a, w_b) )

N = 500000
i = 0
while i<N:

    temp_pos = current_pos  + \
    np.random.normal(0, sd, np.array(start_point).shape )

    temp_obj_value = objective(temp_pos)

    if  temp_obj_value < current_obj_value:

        temp_powell = minimize(
            objective,
            temp_pos,
            method='Powell'
        )
        temp_powell_pos = temp_powell.x
        temp_powell_value = temp_powell.fun
        if temp_powell_value < temp_obj_value:
            temp_obj_value = temp_powell_value
            temp_pos = temp_powell_pos

        current_pos = temp_pos
        current_obj_value = temp_obj_value

    else:
        crunch_step += 1

    if current_obj_value < target:
        break
    if crunch_step % crunch_period == 0:
        crunch_step = 0

    sd = w_cache[crunch_step]

    i += 1
print(current_pos)
print(objective(current_pos))

print('Time: ', ( datetime.datetime.now() ) - start_time)
```

## Discussion

The hybridization of Gaussian Crunching Search (GCS) and Powell's Method for derivative-free optimization is the focus of this research paper. GCS has shown promise in addressing the challenges faced by traditional derivative-free optimization methods. However, it may not always excel in finding the local minimum. On the other hand, traditional methods such as Powell's Method often perform better in this aspect. In our experiments, we observed that GCS exhibits strength in escaping the trap of local minima and approaching the global minima. By combining GCS with specific traditional derivative-free optimization methods, we were able to achieve significant performance improvements while retaining the advantages of each individual method. This hybrid approach offers new possibilities for optimizing complex systems and finding optimal solutions across various applications. Specifically, the Powell-Gaussian Crunching Search (P-GCS) hybrid approach demonstrates its ability to approach better optima that are closer to the global minimum. In contrast, Powell's Method tends to get trapped at local minima. These findings highlight the potential of hybridization techniques to enhance optimization performance and provide valuable insights for future research in this field.

## Conclusion

In conclusion, this research paper introduces a hybrid approach that combines Gaussian Crunching Search (GCS) with Powell's Method for derivative-free optimization. GCS has proven effective in overcoming challenges faced by traditional methods, but it may struggle to find the local minimum. On the other hand, traditional methods may perform better in this aspect but can be trapped at local minima. By combining these two approaches, we have achieved significant performance improvements while retaining the advantages of each method. The hybridization of GCS and Powell's Method opens up new possibilities for optimizing complex systems and finding optimal solutions across various applications. Specifically, the Powell-Gaussian Crunching Search (P-GCS) approach demonstrates its ability to approach optima closer to the global minimum. Future research can further explore and refine this hybrid approach to enhance optimization algorithms and address real-world optimization challenges.

## Reference

[1] Wong, B. (2023) A new derivative-free optimization method: Gaussian crunching search, arXiv.org. Available at: https://arxiv.org/abs/2307.14359 (Accessed: 09 August 2023).
[2] Powell, M.J. (1964) 'An efficient method for finding the minimum of a function of several variables without calculating derivatives', The Computer Journal, 7(2), pp. 155–162. doi:10.1093/comjnl/7.2.155.
[3] Miller, R.E. (2000) Optimization: Foundations and applications. New York: Wiley.